\documentclass[leqno,12pt,twoside]{amsart}
\usepackage{geometry}

\usepackage{times}

\usepackage{amsmath, amssymb, amsfonts, latexsym, mdwlist}
\usepackage[all]{xy}

\usepackage{paralist}
\setdefaultenum{(a)}{(i)}{}{}
\usepackage{enumitem} 

\def\C{\ensuremath{\mathbb{C}}}
\def\D{\ensuremath{\mathbb{D}}}

\def\P{\ensuremath{\mathbb{P}}}
\def\Q{\ensuremath{\mathbb{Q}}}
\def\R{\ensuremath{\mathbb{R}}}

\def\ch{\mathop{\mathrm{ch}}\nolimits}

\def\Coh{\mathop{\mathrm{Coh}}\nolimits}

\def\cok{\mathop{\mathrm{cok}}}

\def\dim{\mathop{\mathrm{dim}}\nolimits}

\def\Ext{\mathop{\mathrm{Ext}}\nolimits}

\def\Hom{\mathop{\mathrm{Hom}}\nolimits}

\def\RlHom{\mathop{\mathbf{R}\mathcal Hom}\nolimits}

\def\NS{\mathop{\mathrm{NS}}\nolimits}

\def\rk{\mathop{\mathrm{rk}}}

\newenvironment{Prf}{\textit{Proof.}\/}{\hfill$\Box$}

\def\MG13{\ensuremath{{\mathcal M}_{\Gamma_1(3)}}}
\def\tildeMG13{\ensuremath{\widetilde{\mathcal M}_{\Gamma_1(3)}}}

\def\into{\ensuremath{\hookrightarrow}}
\def\onto{\ensuremath{\twoheadrightarrow}}

\def\blank{\underline{\hphantom{A}}}

\def\Db{\mathrm{D}^{\mathrm{b}}}

\newcommand\stv[2]{\left\{#1\,\colon\,#2\right\}}

\makeatletter
\newtheorem*{rep@theorem}{\rep@title}
\newcommand{\newreptheorem}[2]{%
\newenvironment{rep#1}[1]{%
 \def\rep@title{#2 \ref{##1}}%
 \begin{rep@theorem}}%
 {\end{rep@theorem}}}
\makeatother

\newtheorem{Thm}{Theorem}[section]
\newreptheorem{Thm}{Theorem}

\newtheorem{Prop}[Thm]{Proposition}
\newtheorem{Lem}[Thm]{Lemma}
\newtheorem{Cor}[Thm]{Corollary}
\newreptheorem{Cor}{Corollary}
\newtheorem{Con}[Thm]{Conjecture}

\theoremstyle{definition}
\newtheorem{Def-s}[Thm]{Definition}
\newtheorem{Def}[Thm]{Definition}
\newtheorem{Rem}[Thm]{Remark}

\def\C{\ensuremath{\mathbb{C}}}
\def\D{\ensuremath{\mathbb{D}}}

\def\P{\ensuremath{\mathbb{P}}}
\def\Q{\ensuremath{\mathbb{Q}}}
\def\R{\ensuremath{\mathbb{R}}}

\def\BB{\ensuremath{\mathcal B}}

\def\FF{\ensuremath{\mathcal F}}

\def\JJ{\ensuremath{\mathcal J}}

\def\OO{\ensuremath{\mathcal O}}

\def\TT{\ensuremath{\mathcal T}}

\def\cht{\ensuremath{\ch^B}}

\newcommand{\ignore}[1]{}

\begin{document}
\title[Stability conditions on threefolds and
Fujita's conjecture]{Bridgeland Stability conditions on threefolds II:
An application to Fujita's conjecture}

\author{Arend Bayer}
\address{Department of Mathematics, University of Connecticut U-3009, 196 Auditorium Road,
Storrs, CT 06269-3009, USA}
\email{bayer@math.uconn.edu}
\urladdr{http://www.math.uconn.edu/~bayer/}

\author{Aaron Bertram}
\address{Department of Mathematics, University of Utah, 155 S 1400 E, Salt Lake City, UT 84112, USA}
\email{bertram@math.utah.edu}
\urladdr{http://www.math.utah.edu/~bertram/}

\author{Emanuele Macr\`i}
\address{Mathematical Institute, University of Bonn, Endenicher Allee 60, D-53115 Bonn, Germany \& Department of Mathematics, University of Utah, 155 S 1400 E, Salt Lake City, UT 84112, USA}
\curraddr{Department of Mathematics, The Ohio State University, 231 W 18th Avenue, Columbus, OH 43210, USA}
\email{macri.6@math.osu.edu}
\urladdr{http://www.math.osu.edu/~macri.6/}

\author{Yukinobu Toda}
\address{Institute for the Physics and Mathematics of the Universe, University of Tokyo, 5-1-5 Kashiwanoha, Kashiwa, 277-8583, Japan}
\email{yukinobu.toda@ipmu.jp}

\keywords{Bogomolov-Gieseker inequality, Bridgeland stability conditions,
Derived category, adjoint line bundles, Fujita conjecture
}

\subjclass[2000]{14F05 (Primary); 14C20, 14J30, 14J32, 18E30 (Secondary)}
\date{\today}

\begin{abstract}
We apply a conjectured inequality 
on third chern classes of stable two-term complexes on threefolds to Fujita's
conjecture. More precisely, the inequality is shown to imply a Reider-type theorem
in dimension three which in turn implies that $K_X + 6L$ is very ample when $L$ is ample, and that
$5L$ is very ample when $K_X$ is trivial.
\end{abstract}

\maketitle

\setcounter{tocdepth}{1}

\tableofcontents

\section{Introduction}\label{sec:intro}

A Bogomolov-Gieseker-type inequality on Chern classes of ``tilt-stable" objects in the derived 
category of a threefold was conjectured in \cite{BMT:3folds-BG} 
in the context of constructing Bridgeland stability conditions. 
In this paper, we show how the same inequality would allow one to
extend Reider's stable-vector bundle technique (\cite{Reider:Surfaces}) from surfaces to threefolds, 
and in particular to obtain Fujita's conjecture in the 
threefold case. This follows a line of reasoning that was suggested in \cite{AB:Reider}. 

While we use the setup of tilt-stability from \cite{BMT:3folds-BG}, this
paper is intended to be self-contained, and to be readable by birational geometers with a
passing familiarity with derived categories.

Tilt-stability depends on two numerical parameters: an ample class $\omega \in \NS_\Q(X)$ and
an arbitrary class $B \in \NS_\Q(X)$. It is a notion of stability on a particular abelian category,
$\BB_{\omega,B}$, of two-term complexes in $\Db(X)$, and codimension three Chern classes of 
stable objects $E$ in this category 
(and not stable vector bundles) are conjectured to satisfy a Bogomolov-Gieseker inequality in Conjecture  \ref{con:strong-BG}. 
Assuming this conjecture, we prove the following Reider-type theorem for threefolds:

\begin{repThm}{thm:final}
Let $X$ be a smooth projective threefold over $\C$, and let $L$ be an ample line bundle on $X$ such that Conjecture \ref{con:strong-BG} holds when $B$ and $\omega$ are scalar multiples of
$L$.
Fix a positive integer $\alpha$, and assume that $L$ satisfies the following conditions:
\begin{enumerate}[label={(\Alph*)}]
\item \label{enum:vol} $L^3> 49\alpha$;
\item \label{enum:divisors} $L^2.D\geq7\alpha$, for all integral divisor classes $D$ with $L^2.D > 0$ and $L.D^2<\alpha$;
\item \label{enum:curves} $L.C\geq3\alpha$, for all curves $C$.
\end{enumerate}
Then $H^1(X,K_X \otimes L \otimes I_Z)=0$ for any zero-dimensional subscheme $Z \subset X$ of length $\alpha$.
\end{repThm}

Theorem \ref{thm:final} would give an effective numerical criterion for an adjoint line bundle to be globally generated ($\alpha = 1$) or very ample ($\alpha = 2$):

\begin{Cor}[Fujita's Conjecture] \label{thm:Fujita}
Let $L$ be an ample line bundle on a smooth projective threefold $X$.
Assume Conjecture \ref{con:strong-BG} holds for $\omega$ and $B$ as above.
Then:
\begin{enumerate}
\item $K_X \otimes L^{\otimes m}$ is globally generated for $m \ge 4$. Moreover, if $L^3\geq2$, then
$K_X \otimes L^{\otimes 3}$ is also globally generated. \label{enum:Fujita1}
\item $K_X \otimes L^{\otimes m}$ is very ample for $m \ge 6$. \label{enum:Fujita2}
\end{enumerate}
\end{Cor}

In Proposition \ref{prop:FujitaL5}, we also show (assuming the conjecture) that $K_X \otimes L^5$ is very ample
as long as its restriction to special degree one curves is very ample. As a consequence, $K_X \otimes L^5$
is very ample when $K_X$ is trivial, or, more generally, when $K_X . C$ is even for all curves $C \subset X$. 

Ein and Lazarsfeld proved that $K_X \otimes L^{\otimes 4}$ is globally generated \cite{EinLazarsfeld:Fujita}.
In the case $L^3\geq 2$, Fujita, Kawamata, and Helmke proved that $K_X \otimes L^{\otimes 3}$ is globally generated as well \cite{Fujita:Threefolds, Kawamata:FujitaConj, Helmke:Fujita}.
In fact, in Proposition \ref{prop:FujitaConverse}, we show that these results conversely give some evidence for Conjecture \ref{con:strong-BG}.
Case \eqref{enum:Fujita2} in Corollary \ref{thm:Fujita} instead is not known in general; but also
note that the strongest form of Fujita's conjecture predicts that $K_X \otimes L^{\otimes 5}$ is already very ample.
For further references, we refer to \cite[Section 10.4]{Laz:Positivity2}.
Notice that the bounds in Theorem \ref{thm:final} are very similar to those in \cite{Fujita:Threefolds} when $\alpha=1$ (see
also \cite{Kawamata:FujitaConj, Helmke:Fujita}) and, when $\alpha=2$ and $Z$ consists  of two
distinct points, to those in \cite{Fujita:VeryAmple}.

\subsection*{Approach}
We explain our approach, which was outlined in \cite[Section 5]{AB:Reider}, but can now be made precise using the strong Bogomolov-Gieseker conjecture of \cite{BMT:3folds-BG}.
It is closer to Reider's original approach \cite{Reider:Surfaces} for surfaces via stability of
sheaves (generalized to threefolds by extending it to derived categories), than to the
Ein-Lazarsfeld-Kawamata approach mentioned above, via vanishing theorems.

Let us give first a brief recall on Reider's method for proving Fujita's Conjecture in the case of $X$ being a surface.
By Serre duality, an adjoint linear system $K_X \otimes L$ is very ample if and only if $\Ext^1(L
\otimes I_Z, \OO_X)=H^1(X, K_X \otimes L \otimes I_Z)^\vee = 0$, for all zero-dimensional subscheme
$Z \subset X$ of length one or two.
If this group was non-zero, we would get a rank $2$ torsion-free sheaf $E$ as the non-trivial
extension $\OO_X \into E \onto L\otimes I_Z$.
Reider's idea is to consider the slope-stability of $E$.
If $E$ is stable, then the classical Bogomolov-Gieseker inequality gives a bound on the degree $L^2$ of $L$
in terms of the length of $Z$.
If $E$ is not stable, then the destabilizing subsheaf gives a curve of bounded degree with respect
to $L$.
Hence, if we assume that $L$ satisfies inequalities similar to \ref{enum:vol} and \ref{enum:curves}, we would
get a contradiction.

We generalize this approach to threefolds as follows.
We suppose the conclusion of Theorem \ref{thm:final} is false. Then by Serre duality, 
$$0 \neq \Ext^2(L \otimes I_Z, \OO_X) = \Ext^1(L \otimes I_Z, \OO_X[1]).$$
For appropriate choices of $\omega$ and $B$, both $L \otimes I_Z$ and $\OO_X[1]$ are objects
in the abelian category $\BB_{\omega, B}$, and thus this extension class
corresponds to another object $E$ of $\BB_{\omega, B}$.
In Section \ref{subsec:ReiderBG}, we will show that for $\omega \to 0$, the complex $E$ violates the
inequality of Conjecture \ref{con:strong-BG}, thus it must become unstable.
We show in Section \ref{subsec:ReiderHodge} that the Chern classes
of a destabilizing subobject give a contradiction to Assumptions \ref{enum:vol} and \ref{enum:divisors}
of the Theorem unless it is of the form $L \otimes I_C$, where $I_C$ is the ideal sheaf of a curve containing $Z$.
In Section \ref{sec:Fujita}, we apply our conjecture and Assumption \ref{enum:curves} to this
remaining case and deduce Theorem \ref{thm:final}.

\subsection*{Acknowledgements}
A.~Ba.~ is partially supported by NSF grant DMS-0801356/DMS-1001056.
A.~Be.~ is partially supported by NSF grant DMS-0901128.
E.~M.~ is partially supported by NSF grant DMS-1001482/DMS-1160466, Hausdorff Center for Mathematics, Bonn, and by SFB/TR 45.
Y.~T.~is supported by World Premier International Research Center Initiative (WPI initiative), MEXT, Japan, and Grant-in AId for Scientific Research grant (22684002), partly (S-19104002), from the Ministry of Education, Culture, Sports, Science and Technology, Japan.
The authors would like to thank the Isaac Newton Institute and its program on ``Moduli Spaces'', during which this paper was finished.

\subsection*{Notation and Convention}
Throughout the paper, $X$ will be a smooth projective threefold defined over $\C$ and $\Db(X)$ its bounded derived category of coherent sheaves.
Given a line bundle $L$ on $X$, we will denote by $\D_L \colon \Db(X) \to \Db(X)$ the following
local dualizing functor on its derived category:
\[
\D_L(\blank) := (\blank)^\vee[1] \otimes L = \RlHom(\blank, L[1]).
\]
We identify a line bundle $L$ with its first Chern class $c_1(L)$, and write $K_X$ for the canonical line bundle.
While $L^{\otimes m}$ denotes the tensor powers of the line bundle, $L^k$ denotes the intersection
product of its first Chern class.

\section{Setup}\label{sec:recall}

In this section, we briefly recall the notion of ``tilt-stability'' defined in
\cite[Section 3]{BMT:3folds-BG} and its most important properties.

Let $X$ be a smooth projective threefold, and let $\omega, B \in \NS_{\Q}(X)$ be rational numerical divisor classes such that $\omega$ is ample.
We use $\omega, B$ to define a slope function $\mu_{\omega, B}$ for coherent sheaves on $X$ as follows:
For torsion sheaves $E$, we set $\mu_{\omega, B}(E) = + \infty$, otherwise
\[
\mu_{\omega, B}(E) = \frac{\omega^2 \cht_1(E)}{\omega^3 \cht_0(E)}
 = \frac{\omega^2 \ch_1(E)}{\omega^3 \cht_0(E)} - \frac{\omega^2 B}{\omega^3}
\]
where $\cht(E) = e^{-B} \ch(E)$ denotes the Chern character twisted by $B$ (explicitly, $\cht_0=\rk$, $\cht_1=\mathrm{c}_1-B\rk$, etc.).

A coherent sheaf $E$ is slope-(semi)stable (or $\mu_{\omega,B}$-(semi)stable) if, for all subsheaves $F \into E$, we have
\[
\mu_{\omega, B}(F) < (\le) \mu_{\omega, B}(E/F).
\]
Due to the existence of Harder-Narasimhan filtrations (HN-filtrations, for short) with respect to slope-stability,
there exists a ``torsion pair'' $(\TT_{\omega, B}, \FF_{\omega, B})$ defined 
as follows:
\begin{align*}
\TT_{\omega, B} &= \stv{E \in \Coh X}
{\text{any quotient $E \onto G$ satisfies $\mu_{\omega, B}(G) > 0$}} \\
\FF_{\omega, B} &= \stv{E \in \Coh X}
{\text{any subsheaf $F \into E$ satisfies $\mu_{\omega, B}(F) \le 0$}}
\end{align*}
Equivalently, $\TT_{\omega, B}$ and $\FF_{\omega, B}$ are the extension-closed
subcategories of $\Coh X$ generated by slope-stable sheaves of positive or non-positive slope,
respectively.

\begin{Def}\label{def:BB}
We let $\BB_{\omega, B} \subset \Db(X)$ be the extension-closure
\[ \BB_{\omega, B} = \langle \TT_{\omega, B}, \FF_{\omega, B}[1] \rangle.
\]
\end{Def}
More explicitly, $\BB_{\omega, B}$ is the subcategory of two-term complexes
$E\colon E^{-1} \xrightarrow{d} E^0$ with $H^{-1}(E) = \ker d \in \FF_{\omega, B}$ and
$H^0(E) = \cok d \in \TT_{\omega, B}$.
We can characterize isomorphism classes of objects in $\BB_{\omega,B}$ by extension classes:
to give an object $E\in\BB_{\omega,B}$ is equivalent to giving $T\in\TT_{\omega, B}$, $F\in\FF_{\omega, B}$, and a class $\xi\in\Ext^2_X(T,F)$.

By the general theory of torsion pairs and tilting \cite{Happel-al:tilting},
$\BB_{\omega, B}$ is the heart of a bounded t-structure on $\Db(X)$.
For the most part, we only need that $\BB_{\omega, B}$ is an abelian category: Exact sequences in
$\BB_{\omega, B}$ are given by exact triangles in $\Db(X)$.  For any such exact sequence
\[
0\to E\to F\to G\to 0
\]
in $\BB_{\omega,B}$, we have a long exact sequence in $\Coh X$:
\begin{alignat*}{4}
0 &\to H^{-1}(E) & &\to H^{-1}(F) & &\to H^{-1}(G) & &\to \\
  &\to H^{0}(E)  & &\to H^{0}(F)  & &\to H^{0}(G)  & &\to 0.
\end{alignat*}


Using the classical Bogomolov-Gieseker inequality and Hodge Index theorem, we
defined the following slope function on $\BB_{\omega, B}$:
We set $\nu_{\omega, B}(E) = +\infty$ when $\omega^2 \cht_1(E) = 0$, and
otherwise
\begin{equation} \label{eq:nu-def}
\nu_{\omega, B}(E) =
\frac{\omega \cht_2(E) - \frac 16 \omega^3 \cht_0(E)}{\omega^2 \cht_1(E)}.
\end{equation}
We showed that this is a slope function, in the sense that it satisfies the weak see-saw property
for short exact sequences in $\BB_{\omega, B}$: for any subobject $F \into E$, we have
$\nu_{\omega, B}(F) \le \nu_{\omega, B}(E) \le \nu_{\omega, B}(E/F)$ or
$\nu_{\omega, B}(F) \ge \nu_{\omega, B}(E) \ge \nu_{\omega, B}(E/F)$.

\begin{Def} \label{def:tilt-stable}
An object $E \in \BB_{\omega, B}$ is ``tilt-(semi)stable'' if, for all non-trivial subobjects
$F \into E$, we have
\[
\nu_{\omega, B}(F) < (\le) \nu_{\omega, B}(E/F).
\]
\end{Def}

Motivated by the case of torsion sheaves (\cite[Proposition 7.1.1]{BMT:3folds-BG}), by projectively flat vector bundles (\cite[Proposition 7.4.2]{BMT:3folds-BG}), and
the case of $X = \P^3$ (\cite[Theorem 8.2.1]{BMT:3folds-BG} and \cite{Macri:P3}), we stated the following conjecture:

\begin{Con}[{\cite[Conjecture 1.3.1]{BMT:3folds-BG}}] \label{con:strong-BG}
For any $\nu_{\omega, B}$-semistable object $E\in \BB_{\omega, B}$ satisfying 
$\nu_{\omega, B}(E) = 0$, we have the following inequality
\begin{align} \label{eq:strong-BG}
\cht_3(E) \le \frac{\omega^2}{18}\cht_1(E). 
\end{align}
\end{Con}

Conjecture \ref{con:strong-BG} is analogous to the classical Bogomolov-Gieseker inequality, which 
can be formulated as follows:
For any $\mu_{\omega,B}$-semistable sheaf $E$ satisfying $\mu_{\omega,B}(E)=0$, we have $\omega\cht_2(E)\le 0$.

The original motivation for Conjecture \ref{con:strong-BG} is to construct examples of Bridgeland stability conditions on $\Db(X)$.
While any linear inequality of the form \eqref{eq:strong-BG} would be sufficient to this end,
the constant $\frac 1{18}$ in equation \eqref{eq:strong-BG} is chosen so that, if $\omega$ and $B$ are proportional to the first Chern class
of an ample line bundle $L$, the inequality is an equality for tensor power
$L^{\otimes n}$ of $L$. More generally, it is an equality when
$E$ is a slope-stable vector bundles $E$ whose discriminant $\Delta=(\cht_1)^2-2\cht_0\cht_2$ satisfies
$\omega \Delta(E) = 0$, and for which $\cht_1(E)$ is proportional to $L$. Such vector bundles 
have a projectively flat connection, and are examples of tilt-stable objects:

\begin{Prop}[{\cite[Proposition 7.4.1]{BMT:3folds-BG}}] \label{prop:StabilityLineBundles}
Let $L$ be an ample line bundle, and assume that both $\omega$ and $B$ are proportional to $L$. 
Then any slope-stable vector bundle $E$, with
$\omega \Delta(E) = 0$ and for which $\cht_1(E)$ is proportional to $L$, is also tilt-stable
with respect to $\nu_{\omega, B}$.
\end{Prop}

The proof is essentially the same as for line bundles $L^{\otimes n}$ in
\cite[Proposition 3.6]{AB:Reider}. 

By assuming Conjecture \ref{con:strong-BG}, we can also show conversely: if an object in $\BB_{\omega,B}$ is tilt-stable and the inequality in Conjecture \ref{con:strong-BG} is an equality, then it must have trivial discriminant.
We first recall that, based on Bridgeland's deformation theorem in \cite{Bridgeland:Stab}, we also showed
the existence of a continuous family of stability conditions depending on \emph{real} classes $\omega, B$:

\begin{Prop}[{\cite[Corollary 3.3.3]{BMT:3folds-BG}}] \label{prop:openness}
Let $U \subset \NS_\R(X) \times \NS_\R(X)$ be the subset of pairs of real 
classes $(\omega, B)$ for which $\omega$ is ample.
There exists a notion of ``tilt-stability'' for every $(\omega, B) \in U$.
For every object $E$, the set of $(\omega, B)$ for which $E$ is
$\nu_{\omega, B}$-stable defines an open subset of $U$.
\end{Prop}

By using Proposition \ref{prop:openness}, we can then prove the following.

\begin{Prop}\label{prop:ConjEqual}
Let $L$ be an ample line bundle, and assume that both $\omega$ and $B$ are proportional to $L$.
Assume also that Conjecture \ref{con:strong-BG} holds for such $B$ and $\omega$.
Let $E\in\BB_{\omega,B}$ be a $\nu_{\omega, B}$-stable object, with $\ch_0(E)\neq0$ and $\ch_1(E)$
proportional to $L$, and satisfying:
\[
\frac{\omega^3}{6} \ch_0(E)=\omega \ch_2^B(E)\qquad \text{and}\qquad
\ch_3^B(E)=\frac{\omega^2}{18}\ch_1^B(E).
\]
Then $\omega.\Delta(E)=0$.
\end{Prop}

\begin{Prf}
Write $d=L^3$,  $B=b_0L$, $\omega=T_0 L$ and  $\ch_0(E)=r$.
The idea for the proof is that, since stability is an open property, we can deform 
$b=b_0$ and $T = T_0$, as a function $T = T(b)$ of $b$, slightly such that $E$ is still 
$\nu_{T(b)L, bL}$-stable with $\nu_{T(b)L, bL}(E) = 0$; then we apply Conjecture \ref{con:strong-BG}
for the pairs $\omega = T(b)L, B = bL$ depending on $b$.

Evidently, $\nu_{TL, bL}(E) = 0$ is equivalent to
\[
T^2=\frac 6{rd} L.\ch_2^{bL}(E)
\]
Since $T_0 > 0$, and since the equation is satisfied for $T = T_0$ and $b = b_0$, the equation
defines a function $T = T(b)$ for $b$ nearby $b_0$.

It is immediate to check from the definition that the chain rule
\begin{equation}		\label{eq:neatderivative}
\frac{\partial}{\partial b} \ch_i^{bL}(E) = -L\ch_{i-1}^{bL}(E)
\end{equation}
holds for $i = 1, \dots, 3$.

Consider 
\begin{equation*}
f(b) = \ch_3^{bL}(E) - \frac{(T(b)L)^2}{18}.\ch_1^{bL}(E)
= \ch_3^{bL}(E) - \frac{1}{3rd}L.\ch_2^{bL}(E) \cdot L^2.\ch_1^{bL}(E)
\end{equation*}
as a function of $b$ in some neighborhood of $b_0 \in \R$.
By Proposition \ref{prop:openness} and Conjecture \ref{con:strong-BG}, we have $f(b) \le 0$ for $b$
close to $b_0$, and by assumption $f(b_0) = 0$; therefore $f'(b_0) = 0$. Using equation
\eqref{eq:neatderivative}, we obtain
\begin{equation*}
\begin{split}
f'(b) &= -L.\ch_2^{bL}(E) + \frac 1{3rd}\bigl((L^2.\ch_1^{bL})^2 + L.\ch_2^{bL}(E) \cdot rd\bigr)
\\
&= \frac 1{3r} \bigl(L.(\ch_1^{bL}(E))^2 - 2 L.\ch_2^{bL}(E) r\bigr) = \frac 1{3r} L.\Delta(E).
\end{split} \end{equation*}
(Note that we used $(L^2.\ch_1^{bL})^2 = L^3 \cdot L.(\ch_1^{bL})^2$, which
holds because $\ch_1^{bL}(E)$ is proportional to $L$.) This proves the claim.
\end{Prf}

Finally, based on an alternate construction of tilt-stability, we also showed that it behaves
well with respect to the dualizing functor $\D_L(\blank) = \RlHom(\blank, L[1])$ for every
line bundle $L$. For this purpose, we fix $B = \frac L2$:

\begin{Prop}
\label{prop:selfdual}
Let $F \in \BB_{\omega, \frac L2}$ be an object with $\nu_{\omega, \frac L2}(A) < +\infty$
for every subobject $A \subset F$. Then there is an exact triangle
$\widetilde F \to \D_L(F) \to T_0[-1]$
where $T_0$ is a zero-dimensional torsion sheaf and $\widetilde F$ an object of
$\BB_{\omega, \frac L2}$ with
$\nu_{\omega, \frac L2}(\widetilde F) = - \nu_{\omega, \frac L2}(F)$. The object
$\widetilde F$ is $\nu_{\omega, \frac L2}$-semistable if and only if $F$ is $\nu_{\omega, \frac
L2}$-semistable.
\end{Prop}
\begin{Prf}
Since $\D_L(\blank)$ can be written as the composition $\blank \otimes L \circ \D(\blank)$,
this follows from \cite[Proposition 5.1.3]{BMT:3folds-BG} and the fact that tensoring with
$L$ corresponds to replacing $B$ with $B-L$.
\end{Prf}

\section{Reduction to curves}

In this section, we use Assumptions \ref{enum:vol} and \ref{enum:divisors} of Theorem \ref{thm:final}
to show that the non-vanishing of $H^1(X,K_X\otimes L\otimes I_Z)$ implies the existence
of special low-degree curves on $X$.
The approach, explained in the introduction, involves studying the tilt-stability of a certain object $E$ in the category $\BB$ constructed in the previous section.

\subsection{Bogomolov-Gieseker inequalities and stability}\label{subsec:ReiderBG}
We will use Conjecture \ref{con:strong-BG} in the case where $L$ is an ample line bundle on $X$,
$\omega = TL$ for some $T > 0$, and $B = \frac L2$.
The abelian category $\BB:=\BB_{TL, \frac L2}$ is independent of $T$.

To simplify notation, we will rescale the slope function: set $t = \frac{T^2}{6}$ and write $\nu_t$ for
\begin{equation} \label{eq:nut-def}
\nu_t(\blank) = T \cdot \nu_{T L, \frac L2}(\blank)
= \frac{L.\ch^{L/2}_2(\blank) - t d \ch^{L/2}_0(\blank)} {L^2 . \ch^{L/2}_1(\blank)},
\end{equation}
where $d:=L^3$.
Then the inequality of Conjecture \ref{con:strong-BG} states that, for every
$\nu_t$-stable object $E$, we have 
\begin{equation} \label{eq:strong-BGt}
\ch^{L/2}_3(E) \le \frac{t}{3} L^2 . \ch^{L/2}_1(E) \quad \text{if} \quad
L . \ch^{L/2}_2(E) = d t \ch^{L/2}_0(E).
\end{equation}

Let $Z\subset X$ be a zero-dimensional subscheme of length $\alpha$.
Following \cite{AB:Reider}, observe that if $H^1(X, K_X \otimes L \otimes I_Z) \neq 0$, then 
by Serre duality, we also have $\Ext^2(L \otimes I_Z, \OO_X) \neq 0$.
Any non-zero element
$\xi \in \Ext^2(L \otimes I_Z, \OO_X)$ gives a non-trivial exact triangle in $\Db(X)$
\begin{equation}\label{eq:extension}
\OO_X[1]\to E=E_{\xi} \to L\otimes I_Z\xrightarrow{\xi}\OO_X[2].
\end{equation}

We will show that $E$ is $\nu_{t}$-semistable for $t = \frac 18$;
its Chern classes invalidate the inequality of Conjecture 
\ref{con:strong-BG} for $t \ll 1$, and thus it must become unstable for $t < t_0$ and some
$t_0 \in (0, \frac 18]$;
finally, we will show that the the Chern classes of its destabilizing factor would give special curves or
divisors on $X$.

\begin{Prop}\label{prop:StabilityOfE}
Assume that $H^1(X, K_X \otimes L \otimes I_Z) \neq 0$, and let $E$ be an extension as given by equation \eqref{eq:extension}.
\begin{enumerate}
\item\label{eq:Stab1}$E\in\BB$ and
\[
\ch^{L/2}(E)=\left(0,L,0,\frac{d}{24}-\alpha \right).
\]
\item \label{enum:unstable18}
If $t > \frac 18$, then \eqref{eq:extension} destabilizes $E$ with respect to $\nu_t$.
\item\label{eq:Stab3} If $t=\frac 18$, then $E$ is $\nu_t$-semistable.
\item \label{enum:small-t-unstable}
Assume Conjecture \ref{con:strong-BG} and Assumption \ref{enum:vol} of Theorem
\ref{thm:final}.
Then $E$ is not $\nu_t$-semistable for $0<t\ll 1$, 
\end{enumerate}
\end{Prop}

\begin{Prf}
First of all, we have
\begin{align*}
\ch^{L/2}(\OO_X) &= \left(1,-\frac{L}{2},\frac{L^2}{8},-\frac{L^3}{48}\right),\\
\ch^{L/2}(L\otimes I_Z) &= \left(1,\frac{L}{2},\frac{L^2}{8},\frac{L^3}{48}-\alpha\right).
\end{align*}
As $\OO_X$ and $L\otimes I_Z$ are slope-stable, with $\mu_{\omega, L/2}(\OO_X) < 0$ and
$\mu_{\omega, L/2}(L \otimes I_Z) > 0$, we have $\OO_X \in \FF$ and $L \otimes I_Z \in \TT$. By
the definition of $\BB$, it follows that
$\OO_X[1]$, $L\otimes I_Z$ and $E$ are all objects of $\BB$; in particular, we have proved \eqref{eq:Stab1}.

Moreover, we have
\begin{equation} \label{eq:nu-OX-L}
\nu_t(\OO_X[1]) = 2\left( t-\frac 18\right), \quad \nu_t(E) = 0
\end{equation}
which immediately implies \eqref{enum:unstable18}, since \eqref{eq:extension} is an exact sequence in $\BB$.

To prove \eqref{eq:Stab3}, simply observe that, by Proposition \ref{prop:StabilityLineBundles}, both $\OO_X[1]$ and $L$ are $\nu_t$-stable for all $t>0$.
Moreover, since $\nu_t(L \otimes I_Z) = \nu_t(L)$, any destabilizing subobject $A \into L \otimes I_Z$
would also destabilize $L$ via the composition
$A \into L \otimes I_Z \into L$ (which is an inclusion in $\BB$); thus
$L \otimes I_Z$ is also $\nu_t$-stable. For $t = \frac 18$, we have
$\nu_t(\OO_X[1]) = \nu_t(L \otimes I_Z) = 0$, and thus the extension \eqref{eq:extension} shows
that $E$ is $\nu_t$-semistable at $t = \frac 18$.

Finally, if $E$ was $\nu_t$-semistable for all
$t\in(0,\frac 18]$, then by our conjectural inequality \eqref{eq:strong-BGt} we would get
\begin{equation}\label{eqn:Bonn5111}
\frac d{24}-\alpha\leq\frac t{3} d
\end{equation}
for all such $t$.
Hence $d\leq 24\alpha$, in contradiction to Assumption \ref{enum:vol}.
\end{Prf}

Notice that the previous proposition would answer Question 4 in \cite{AB:Reider}.
Also observe that in part \eqref{enum:small-t-unstable},
instead of Assumption \ref{enum:vol}, already assuming $d > 24 \alpha$ would have been enough.
Similarly, instead of Conjecture \ref{con:strong-BG}, any linear inequality between $\cht_3$ and $\cht_1$ would have been sufficient.

In the following proposition, we will show that our situation is self-dual with respect to the local dualizing functor $\D_L(\blank) = \RlHom(\blank, L[1])$.
As a preliminary, let us first note that we may make the following assumption:
\begin{equation*}\tag{*}\label{assptn:star}
H^1(X, K_X\otimes L \otimes I_{Z'}) = 0 \text{ for all subschemes } Z' \subsetneq Z, \text {and } H^1(X, K_X \otimes L \otimes I_Z) \cong \C. 
\end{equation*}
Indeed, in order to show $H^1(X, L \otimes I_Z \otimes  K_X)  = 0$, we can proceed by induction on the length of $Z$ (the case $\alpha = 0$ is, of course, given by Kodaira vanishing).

\begin{Prop}  	\label{prop:E-selfdual}
If Assumption (*) holds, and $E$ is given by the unique non-trivial extension of the form
\eqref{eq:extension}, then $E \cong \D_L(E)$.
\end{Prop}
\begin{Prf}
Due to Assumption (*), it is sufficient to show that $\D_L(E)$ is again a non-trivial extension of
the form \eqref{eq:extension}.
Applying the octahedral axiom to the composition $\OO_Z[-1] \to L \otimes I_Z \to \OO_X[2]$, and
using the two exact triangles \eqref{eq:extension} and
$O_Z[-1] \to L \otimes I_Z \to L$,
we obtain an exact triangle $F \to E \to L$, where $F$ itself fits into
an exact triangle 
\begin{equation} \label{eq:triangle-A}
\OO_X[1] \to F \to \OO_Z[-1].
\end{equation} 

We claim that 
$\Hom(k(x)[-1], F) = 0$ for all skyscraper sheaves of points $x \in X$. Using the long exact
sequence for $\Hom(k(x), \blank)$ applied to \eqref{eq:triangle-A},
we see that this is equivalent to the non-vanishing of
the composition
\begin{equation} \label{eq:inclusion}
k(x)[-1] \to \OO_Z[-1] \to L \otimes I_Z \xrightarrow{\xi} \OO_X[2]
\end{equation} for every inclusion
$k(x) \into \OO_Z$. Given such an inclusion, let $Z' \subset Z$ be the subscheme
given by $\OO_{Z'} \cong \OO_Z/k(x)$. 
If the composition
\eqref{eq:inclusion} vanishes, then $\xi$ factors via $L \otimes I_{Z} \into L \otimes I_{Z'}$.
This contradicts our assumption 
$\Ext^2(L \otimes I_{Z'}, \OO_X) = H^1(X, L \otimes I_{Z'} \otimes  K_X)^\vee = 0$.

Now we apply $\D_L$ to the exact triangle $\OO_X[1] \to F \to \OO_Z[-1]$.
As $\D_L(\OO_X[1]) = L$ and
$\D_L(\OO_Z[-1]) = \OO_Z[-1]$, dualizing \eqref{eq:triangle-A} gives an  exact triangle 
$\OO_Z[-1] \to \D_L(F) \to L \to \OO_Z$. Since
$\Hom(\D_L(F), k(x)[-1]) = \Hom(k(x)[-1], F) = 0$ for all $x \in X$,
the map $L \to \OO_Z$ must be surjective, and hence $\D_L(F) \cong L \otimes I_Z$.
Consequently, applying $\D_L$ to the exact triangle $F \to E \to L$ shows that
$\D_L(E)$ is indeed a non-trivial extension of the form \eqref{eq:extension}.
\end{Prf}

\subsection{Chern classes of destabilizing subobjects}\label{subsec:ReiderHodge}

By Proposition \ref{prop:StabilityOfE} and Proposition \ref{prop:openness}, Conjecture \ref{con:strong-BG} implies the existence of $t_0\in(0,\frac 18]$ with the
following properties:
\begin{itemize}
\item $E$ is $\nu_{t_0}$-semistable.
\item There exists an exact sequence in $\BB$
\begin{equation} \label{eq:destab-E}
0\to A\to E\to F\to 0,
\end{equation}
with $\nu_{t}(A)>0$ if $t<t_0$, and $\nu_{t_0}(A)=0$.
\end{itemize}

In the remainder of this section, we will prove the following statement:

\begin{Prop} \label{prop:reduce-to-curves}
Assume that  $X, L, \alpha$ satisfy 
Assumptions \ref{enum:vol} and \ref{enum:divisors} of Theorem \ref{thm:final} and Assumption (*) 
of the previous section. Then in any destabilizing sequence \eqref{eq:destab-E}, the object
$A$ is of the form
$L \otimes I_C$, for some purely one-dimensional subscheme  $C \subset X$ containing $Z$.
\end{Prop}

We will first prove this for subobjects satsfying
$L^2 . \ch^{L/2}_1(A) \le L^2 . \ch^{L/2}_1(F)$, or, equivalently,
\begin{equation} 		\label{eq:cht1A}
L^2 . \ch^{L/2}_1(A) \le \frac 12 L^2 . \ch^{L/2}_1(E) = \frac d2.
\end{equation}
(We will later use the derived duality $\D_L(\blank)$ to reduce to this case.)

\begin{Lem} \label{lem:rank}
Any subobject $A$ satisfying \eqref{eq:cht1A} is a sheaf with $\rk(A) = \rk(H^0(A)) > 0$.
\end{Lem}
\begin{Prf}
Consider the long exact cohomology sequence for $A \into E \onto F$.
If $H^{-1}(A) \neq 0$, then $H^{-1}(A)\into\OO_X$ is isomorphic to an ideal sheaf of some subscheme $Y$ of $X$.
Since $\OO_Y\into H^{-1}(F)$ and $H^{-1}(F)$ is torsion-free, we must have $H^{-1}(A)\cong\OO_X$.
Then $H^0(A)$ is also torsion-free,
and \eqref{eq:cht1A} implies
\[
L^2 . \ch^{L/2}_1(H^0(A)) = L^2 . \ch^{L/2}_1(A) - L^2 . \ch^{L/2}_1(\OO_X[1]) \le \frac d2 - \frac d2 = 0.
\]
On the other hand, by construction of $\BB$, every HN-filtration factor $U$ of $H^0(A)$
satisfies $L^2 . \ch^{L/2}_1(U) > 0$; thus $H^0(A) = 0$ and $A = \OO_X[1]$. This contradiction
proves $H^{-1}(A) = 0$.

Finally, note that if $A = H^0(A)$ is a torsion-sheaf, then $\nu_t(A)$
is independent of $t$, again a contradiction.
\end{Prf}

\begin{Lem}\label{lem:TorsionFujita}
Either $A$ is torsion-free, or its torsion-part $A_t$ satisfies
\[
L^2.\ch_1(A_t) - 2L.\ch_2(A_t) \geq 0 \quad \text{and} \quad
L^2. \ch_1(A_t) > 0.
\]
\end{Lem}

\begin{Prf}
The sheaf $A_t$ is a subobject of $E$ in $\BB$ with $\rk=0$.
Hence $L.\ch^{L/2}_2(A_t)\leq0$, otherwise it would destabilize $E$ at $t=\frac 18$.
Expanding $\ch^{L/2}_2$ gives the first inequality. To show the second inequality, we just observe that
there are no non-trivial morphisms from sheaves supported in dimension
$\le 1$ to $E$.
\end{Prf}

\begin{Lem} \label{lem:I-II}
In the HN-filtration of $A$ with respect to slope-stability, there exists a factor $U$ of rank $r$ such that
$\Gamma:= L - \frac{\ch_1(U)}r$ satisfies the following inequalities:
\begin{align*}
\label{eq:FundIG}\tag{I} & L^2.\Gamma\leq L.\Gamma^2+6\alpha\\
\label{eq:FundIIG}\tag{II} & \frac d2\left( 1-\frac 1r\right)\leq L^2.\Gamma < \frac d2.
\end{align*}
The case $r = 1$ and $L^2 . \Gamma = 0$ only occurs when $A$ is a torsion-free sheaf
of rank one and $H^{-1}(F) = \OO_X$.
\end{Lem}

If $A$ was a line bundle, the above definition of $\Gamma$ would be just as Reider's original argument
for surfaces: in this case, $\Gamma$ is the support of the cokernel of $A \into H^0(E) \cong L \otimes I_Z$.

\begin{Prf}
From $\nu_{t_0}(A) = 0$ we obtain
\begin{equation} \label{eq:t0}
t_0=\frac{L.\ch^{L/2}_2(A)}{\rk(A)d}.
\end{equation}
Applying the conjectured inequality \eqref{eq:strong-BGt} to $E$, and plugging in $t_0$ gives
\begin{equation*}
\frac d{24}-\alpha= \ch^{L/2}_3(E)\leq\frac{L^2.\ch^{L/2}_1(E)}{3}t_0=\frac d3\frac{L.\ch^{L/2}_2(A)}{\rk(A)d}=\frac 13 \frac{L.\ch^{L/2}_2(A)}{\rk(A)}.
\end{equation*}

We want to bound $L.\ch^{L/2}_2(A)$.  First we expand $\ch^{L/2}_2(A)$:
\[
\ch^{L/2}_2(A)=\ch_2(A)-\frac{L.\ch_1(A)}2+\rk(A)\frac{L^2}8.
\]
Substituting, we deduce
\begin{equation}\label{eq:ineqBonn23sera}
\frac{L^2.\ch_1(A)}{\rk(A)}-2\frac{L.\ch_2(A)}{\rk(A)} \le 6 \alpha.
\end{equation}

Let $A_{tf}$ denote the torsion-free part of $A$, and 
consider its HN-filtration.
Among the HN factors, we choose a torsion-free sheaf $U$ for which the function
\[
\eta(\blank):=\frac{L^2.\ch_1(\blank)-2L.\ch_2(\blank)}{\rk(\blank)}
\]
is minimal.
Notice that $\eta$ satisfies the see-saw property: for an exact sequence of torsion-free sheaves
\[
0\to M\to N\to P\to 0,
\]
we have $\eta(N)\geq\mathrm{min}\{\eta(M),\eta(P)\}$.
Hence we get a chain of inequalities leading to
\begin{equation}\label{eq:eta}
\eta(U)\leq\eta(A_{tf})\leq\eta(A) \le 6\alpha
\end{equation}
where we used Lemma \ref{lem:TorsionFujita} for the second inequality.

To abbreviate, we now write $D:=\ch_1(U)$ and $r:=\rk(U)$.
Since $U$ is $\mu_{\omega, L/2}$-semistable, we can combine the classical Bogomolov-Gieseker inequality 
with \eqref{eq:eta} to obtain
\[
L^2.\frac Dr = \frac{2L.\ch_2(U)}r + \eta(U)\leq L.\frac{D^2}{r^2}+6\alpha.
\]
Substituting $D = rL - r\Gamma$ yields the inequality \eqref{eq:FundIG}.

To prove the chain of inequalities \eqref{eq:FundIIG}, we observe on the one
hand that $L^2.\ch^{L/2}_1(U) > 0$ by the definition of $\TT_{\omega, B} = \BB \cap \Coh X$.
On the other hand, $U$ is a subquotient of $A$ in $\TT_{\omega, B}$; combined
with inequality \eqref{eq:cht1A} we obtain
\[
0 < L^2. \ch^{L/2}_1(U) \leq L^2.\ch^{L/2}_1(A) \leq\frac d2.
\]
Plugging in $\ch^{L/2}_1(U) = -\frac r2 L + D = \frac r2 L - r\Gamma$ shows 
the inequality \eqref{eq:FundIIG}.

Finally, note that in the case $r = 1$ and $L^2. \Gamma = 0$ the chain of inequalities
leading to the first part of \eqref{eq:FundIIG} must be equalities; in particular
$L^2 . \ch^{L/2}_1(U) = L^2 . \ch^{L/2}_1(A)$. This shows that $A_{tf}$ cannot have any other
HN-filtration factors besides $U$, i.e., $U = A_{tf}$. Additionally
it implies that $\ch^{L/2}_1(A_t) = 0$, in contradiction
to Lemma \ref{lem:TorsionFujita}; hence $A_t = 0$ and $A = U$ is a torsion-free rank one sheaf. 

As $L \otimes I_Z$ is torsion-free, if the image of
$H^{-1}(F) \to A$ is non-trival, then the map is surjective, and the inclusion $A \into E$ factors
via $A \into \OO_X[1] \into E$, in contradiction to the stability of $\OO_X[1]$ for all $t$ and
$\nu_{t_0 - \varepsilon}(A) > 0 > \nu_{t_0 - \varepsilon}(\OO_X[1])$. Thus $H^{-1}(F) = \OO_X$.
\end{Prf}

\begin{Prf} (Proposition \ref{prop:reduce-to-curves})
We combine \eqref{eq:FundIG} and \eqref{eq:FundIIG} with the Hodge Index Theorem (just as in
\cite[Corollary 3.9]{AB:Reider}) to obtain
\[
\left(L.\Gamma^2\right)d\leq\left(L^2.\Gamma\right)^2\leq \frac d2\left(L.\Gamma^2+6\alpha\right),
\]
and so $L.\Gamma^2\leq6\alpha$.

In the case $r>1$, we use \eqref{eq:FundIG} and \eqref{eq:FundIIG} again to get
\[
\frac d4\leq L^2.\Gamma\leq L.\Gamma^2 +6\alpha\leq 12\alpha,
\]
and so $d\leq48\alpha$ in contradiction to Assumption \ref{enum:vol}.

Reider's original argument in \cite{Reider:Surfaces} deals with
the case $r=1$:
In case $L^2.\Gamma\neq0$, then $L^2.\Gamma\geq1$.
Let $\kappa:=L.\Gamma^2\leq6\alpha$.
Again combining the Hodge Index Theorem with \eqref{eq:FundIG},
we obtain
\[
\left( L.\Gamma^2\right)d\leq\left( L.\Gamma^2+6\alpha\right)^2,
\]
and so
\[
d\leq12\alpha+\frac{\kappa^2+36\alpha^2}{\kappa}.
\]
The RHS is strictly decreasing function for $\kappa \in (0, 6\alpha]$ and equals
$49 \alpha$ for $\kappa = \alpha$; thus Assumption \ref{enum:vol}
implies $\kappa < \alpha$. On the other hand, $\Gamma$ is integral, and hence Assumption
\ref{enum:divisors} implies $L^2.\Gamma \ge 7\alpha$, in
contradiction to \eqref{eq:FundIG}.

Finally, if $L^2.\Gamma=0$; then, according to Lemma \ref{lem:I-II}, we have $H^{-1}(F)\cong\OO_X$.
Hence $A$ is a subsheaf of $L \otimes I_Z$ with $\ch_1(A) = \ch_1(L)$; this is only possible if
$A \cong L\otimes I_W$, for some closed subscheme $W \subset X$ with $\dim(W)\leq1$.
If $W$ is zero-dimensional, then $\ch^{L/2}_2(A) = \frac 12 L^2$
and equation \eqref{eq:t0} gives $t_0 = \frac 12$, in contradiction to $t_0 \in (0, \frac 18]$.
Hence $W$ is one-dimensional, and
we have shown that any subobject $A$ with $\ch^{L/2}_1(A) \le \frac d2$ is of the form
$A \cong L \otimes I_W$.
In particular $\ch^{L/2}_1(A) = \frac d2$ in this case, so there are no subobject with
$\ch^{L/2}_1(A) < \frac d2$.

Now assume $\ch^{L/2}_1(A) > \frac d2$.  We can apply 
Proposition \ref{prop:E-selfdual} and Proposition \ref{prop:selfdual} to the short exact
sequence \eqref{eq:destab-E} obtain a short exact sequence
in $\BB$
\[
0 \to \widetilde F \xrightarrow{u} E \to E/\widetilde F \to 0
\] 
which is again destabilizing.
Indeed, since $\BB$ is
the heart of a bounded t-structure,
there exists a cohomology functor $H_\BB^*(\blank)$. Applied to the exact triangle
\[
\D_L(F)\to \D_L(E)=E\to \D_L(A),
\]
it induces a long exact sequence in $\BB$
\begin{equation}\label{eq:dual-long}
0\to \widetilde F = H_\BB^0(\D_L(F)) \xrightarrow{u} E \to \widetilde A  \to
T_0 = H_\BB^1(\D_L(F)) \to 0.
\end{equation}

As $\D_L$ preserves $L^2.\ch^{L/2}_1(\blank)$, we have that $\widetilde F$ is a destabilizing
subobject with $\ch^{L/2}_1(F) = \ch^{L/2}_1(E) - \ch^{L/2}_1(A) < \frac d2$, which does not exist.

Finally, note that the long exact sequence \eqref{eq:dual-long} also implies that 
$\D_L(A) = \widetilde A \in \BB$. This gives the vanishing of 
$0 = \Hom(\D_L(A), k(x)[-1]) = \Hom(k(x)[-1], A)$. This is equivalent to the claim that
$W$ is a purely one-dimensional scheme, as any subsheaf $k(x) \into \OO_W$ gives an
extension of $k(x)$  by $L \otimes I_W$.
This finishes the proof of Proposition \ref{prop:reduce-to-curves}.
\end{Prf}

\section{A Reider-type theorem}\label{sec:Fujita} 

In this section we prove our main theorem:

\begin{Thm}\label{thm:final}
Let $L$ be an ample line bundle on a smooth projective threefold $X$, and assume
Conjecture \ref{con:strong-BG} holds for $B$ and $\omega$ proportional to $L$.
Fix a positive integer $\alpha$, and assume that $L$ satisfies
the following conditions:
\begin{enumerate}[label={(\Alph*)}]
\item $L^3> 49\alpha$;
\item $L^2.D\geq7\alpha$, for all integral divisor classes $D$ with
$L^2.D > 0$ and $L.D^2<\alpha$;
\item $L.C \ge 3\alpha$, for all curves $C$.
\end{enumerate}
Then $H^1(X,K_X \otimes L \otimes I_Z)=0$, for any zero-dimensional subscheme $Z \subset X$ of length $\alpha$.
\end{Thm}

\begin{Prf}
As explained in Section \ref{subsec:ReiderBG}, we may proceed by induction on the length of $Z$
and may use Assumption (*).
Let $t_0\in (0,\frac 18]$ be as in Section \ref{subsec:ReiderHodge} and let $t = t_0 - \epsilon$.
Truncating the Harder-Narasimhan filtration of $E$ with respect to
$\nu_t$-stability gives a short exact sequence
\[
0 \to A \to E \to F \to 0
\]
with $\nu_t(A) > 0$, 
such that any subobject $A' \into E$ with $\nu_{t}(A') > 0$ factors via $A' \into A$.
By Proposition \ref{prop:reduce-to-curves}, $A$ is of the form $L \otimes I_C$ for some
purely one-dimensional subscheme $C \subset X$; it also implies that $A$ is stable, as any destabilizing
subobject $A'$ of $A$ would again be of the form $A' \cong L \otimes I_{C'}$, so that the quotient
$A/A'$ would be a torsion sheaf with $\nu_t(A/A') = +\infty$.

Let $\widetilde F$ be the object obtained by dualizing $F$ and applying Proposition
\ref{prop:selfdual}. The map $\D_L(F) \to \D_L(E) \cong E$ induces a map
$\widetilde F \to E$ which is an injection in $\BB$.
Since 
\begin{equation} \label{eq:cht-tilde-F}
\ch^{L/2}_i(\widetilde F) = \ch^{L/2}_i(\D_L(F))
\end{equation}
for $i \le 2$, we have $\nu_t(\widetilde F) = - \nu_t(F) > 0$; thus the map factorizes as
$\widetilde F \into A \into E$. By Proposition \ref{prop:reduce-to-curves}, 
the object $\widetilde F$ is of the form $L \otimes I_{C'}$ for some purely one-dimensional
subscheme $C' \subset X$. Equation \eqref{eq:cht-tilde-F} also implies
$\ch^{L/2}_i(\widetilde F) = \ch^{L/2}_i(A)$ for $i \le 2$; thus the (non-trivial) map
$L \otimes I_{C'} \to L \otimes I_C$ has zero-dimensional cokernel. 
It follows that
\begin{equation*}
\ch^{L/2}_3(F) = \ch^{L/2}_3(\D_L(F)) \le \ch^{L/2}_3(\widetilde F) \le \ch^{L/2}_3(A).
\end{equation*}
This implies that
\begin{equation}\label{eq:ch3bound}
2 \ch^{L/2}_3(A) \ge \ch^{L/2}_3(A) + \ch^{L/2}_3(F) = \ch^{L/2}_3(E) = \frac{d}{24} - \alpha,
\end{equation}
and the difference of the two sides is a non-negative integer.

On the other hand, as $A$ is stable, by Conjecture \ref{con:strong-BG}, by \eqref{eq:t0} and \eqref{eq:ch3bound}, and by expanding $\ch^{L/2}$ we have
\begin{equation}\label{eq:MainCurves}
\frac{d}{48} - \frac{\alpha}{2}\leq\ch^{L/2}_3(A)\leq\frac{t_0}3 L^2.\ch^{L/2}_1(A)=\frac{1}{6}L.\ch^{L/2}_2(A)=\frac{d}{48} - \frac{L.C}6.
\end{equation}

We now use Assumption \ref{enum:curves}: $L.C\geq3\alpha$.
This contradicts \eqref{eq:MainCurves}, unless $L.C=3\alpha$ and
\[
\frac{d}{48} - \frac{\alpha}{2}=\ch^{L/2}_3(A)=\frac{t_0}3 L^2.\ch^{L/2}_1(A).
\]
Since $(TL).\Delta(A) = 3 \alpha T \neq 0$, this in turn
contradicts Proposition \ref{prop:ConjEqual}.
\end{Prf}

We also obtain the following result characterizing the only possible counter-examples to
Fujita's very ampleness conjecture in case $L = M^5$:

\begin{Prop}  \label{prop:FujitaL5}
Assume that Conjecture \ref{con:strong-BG} holds for $X$, $\omega = tL$ and $B = \frac L2$
and $L \cong M^5$ for an ample line bundle $M$.
Then either $ K_X \otimes L$ is very ample, or there exists a curve $C$ of degree
$M.C = 1$ and arithmetic genus $g_a(C) = \frac 52 + \frac 12 K_X.C$ such that $K_X \otimes L|_C$ is a line bundle of degree
$2 g_a(C)$ on $C$ which is not very ample.
\end{Prop}

\begin{Prf}
Assume that $K_X\otimes L$ is not very ample.
We follow the logic and the notation of the proof of Theorem \ref{thm:final}, with $\alpha = 2$. As before,
let $A=L\otimes I_C$ be the
destabilizing subobject of $E$ for $t = t_0 - \epsilon$; here $C$ is a purely one-dimensional subscheme of $X$.
By the proof of Theorem \ref{thm:final}, we have $L.C < 6$ and thus necessarily $M.C = 1$ and $L.C = 5$.
In particular, $C$ is reduced and irreducible.
We claim that $\ch^{L/2}_3(A)=\frac{d}{48} - 1$.
Indeed, setting $\alpha = 2$ in  \eqref{eq:MainCurves} gives
\begin{equation}\label{eq:mfive}
\frac{d}{48} - 1\leq\ch^{L/2}_3(A)\leq\frac{d}{48} - \frac 56.
\end{equation}
On the other hand, if $\ch^{L/2}_3(A)\neq\frac{d}{48} - 1$, then, by \eqref{eq:ch3bound}, $\ch^{L/2}_3(A)\geq\frac{d}{48}-\frac 12$, a contradiction to the inequality \eqref{eq:mfive}.

From the claim, we obtain
\[ \ch_3(L \otimes \OO_C) = \ch_3(L) - \ch_3(A) = \frac 72
\]
and thus
\[ \ch_3(\OO_C) = \ch_3(L \otimes \OO_C) - L.C = - \frac 32
\]
By Hirzebruch-Riemann-Roch, we get
\[
1 - g_a(C) = \ch_3(\OO_C) -\frac 12 K_X.C. 
\]
Plugging in the previous equation and solving for $K_X.C$ shows that $K_X \otimes L|_C$ is a line bundle of degree
$2 g_a(C)$ on $C$.
The explicit expression for $g_a(C)$ follows immediately.

Finally, the cohomology sheaves of the quotient $F\cong E/A$ are $H^{-1}(F) \cong \OO_X$ and $H^0(F)
\cong L \otimes \OO_C(-Z)$ (where $\OO_C(-Z)$ denotes the ideal sheaf of $Z \subset C$).
If $F$ were decomposable, $\widetilde F$ would be a decomposable destabilizing subobject of $E$,
which cannot exist. Hence
\[
0 \neq \Ext^2(L \otimes \OO_C(-Z), \OO_X) = H^1(C, K_X \otimes L|_C(-Z))^\vee.
\]
On the other hand,  $K_X \otimes L|_C$ is a line bundle of degree $2g_a(C)$ on an irreducible
Cohen-Macaulay curve, and thus $H^1(K_X \otimes L|_C) = 0$. Hence
$K_X \otimes L|_C$ is not very ample.
\end{Prf}

\begin{Rem}\label{rmk:FujitaL5}
Notice that Proposition \ref{prop:FujitaL5} implies Fujita's conjecture when $K_X$ is numerically
trivial (or, more generally, when $K_X.C$ is even for all integral curve classes $C$).
\end{Rem}

In case the curve $C \subset X$ of Proposition \ref{prop:FujitaL5} is l.c.i, one can be even more precise.
Let $\omega_C$ be the dualzing sheaf (which agrees with the dualizing complex, as $\OO_C$ is pure
and thus $C$ Cohen-Macaulay). The sheaf $K_X \otimes L(-Z)|_C$ is torsion-free of rank one and degree
$2g_a(C)-2$ with $H^1(K_X \otimes L(-Z)|_C ) \neq 0$, and thus Serre duality implies
$K_X \otimes L(-Z)|_C  \cong \omega_C$. If $N$ is the normal bundle, adjunction gives
$\Lambda^2 N \cong L(-Z)$. In particular, the normal bundle has degree 3. Since $M.C = 1$,
bend-and-break implies that such a curve cannot be rational.

In conclusion, we show how to reverse part of the argument in this section when $Z$ has length one.
Indeed, in such a case we can use Ein-Lazarsfeld theorem (or better, its variant by Kawamata and Helmke) to show that Conjecture \ref{con:strong-BG} holds true for this particular case, coherently with our result:

\begin{Prop}\label{prop:FujitaConverse}
Let $L$ be an ample line bundle on a smooth projective threefold $X$.
Assume that $L$ satisfies the following conditions:
\begin{enumerate}
\item $L^3\geq 28$;
\item $L^2.D\geq9$, for all integral effective divisor classes $D$.
\end{enumerate}
Assume also that there exists $x\in X$ such that $H^1(X,K_X\otimes L\otimes I_{x})\neq0$.
Then Conjecture \ref{con:strong-BG} holds for all objects $E\in\BB$ given as non-trivial extensions
\[
\OO_X[1]\to E \to L\otimes I_x\to \OO_X[2].
\]
\end{Prop}

\begin{Prf}
The argument is very similar to \cite{Kawamata:FujitaConj}, Proposition 2.7 and Theorem 3.1, Step 1.
We freely use the notation from \cite[Sections 9 \& 10]{Laz:Positivity2}.
By \cite[Lemma 2.1]{Kawamata:FujitaConj}, given a rational number $t$ satisfying $3/{\sqrt[3]{L^3}}<t<1$, there exists a $\Q$-divisor $D$ numerically equivalent to $tL$ such that $\mathrm{ord}_x D =3$.
Let $c\leq1$ the log-canonical threshold of $D$ at $x$.

By \cite[Theorem 3.1]{Kawamata:FujitaConj} (also \cite{Helmke:Fujita}) and our assumptions, the LC-locus $\mathrm{LC}(cD;x)$ (i.e., the zero-locus of the multiplier ideal $\JJ(c\cdot D)$ passing through $x$) must be a curve $C$ satisfying $1\leq L.C\leq 2$.
We can now apply Nadel's vanishing theorem to $cD$ to deduce that $H^1(X,K_X\otimes L\otimes I_C)=0$, and so that the restriction map $H^0(X,K_X\otimes L)\onto H^0(X,K_X\otimes L|_C)$ is surjective.

Consider the composition $u\colon L\otimes I_C\to L\otimes I_x\to \OO_x[2]$.
Then, $u\neq 0$ if and only if $x$ is a base point of $K_X\otimes L$ which is not a base point of $K_X\otimes L|_C$.
The surjectivity of the restriction map implies that $u=0$.
Hence, we get an inclusion $L\otimes I_C \into E$ in $\BB$ which destabilizes $E$, if \eqref{eq:strong-BG} is not satisfied.
\end{Prf}


\bibliography{all}                      
\bibliographystyle{alphaspecial}     

\end{document}